\documentclass[a4paper]{article}
\usepackage{amsthm,amsmath,amssymb}
\newtheorem{teor}{Theorem}[section]
\newtheorem{defin}[teor]{Definition}
\newtheorem{lemm}[teor]{Lemma}
\newtheorem{osse}[teor]{Remark}
\newtheorem{prop}[teor]{Proposition}
\newtheorem{defi}[teor]{Definition}
\newtheorem{coro}[teor]{Corollary}
\newtheorem{prob}[teor]{Problem}

\newcommand{\bele}{\begin{lemm}\begin{sl}}
\newcommand{\enle}{\end{sl}\end{lemm}}
\newcommand{\bedef}{\begin{defi}\begin{sl}}
\newcommand{\eddef}{\end{sl}\end{defi}}
\newcommand{\bete}{\begin{teor}\begin{sl}}
\newcommand{\ente}{\end{sl}\end{teor}}
\newcommand{\beos}{\begin{osse}\begin{rm}}
\newcommand{\eddos}{\end{rm}\end{osse}}
\newcommand{\bepr}{\begin{prop}\begin{sl}}
\newcommand{\empr}{\end{sl}\end{prop}}
\newcommand{\bepro}{\begin{prob}\begin{rm}}
\newcommand{\empro}{\end{rm}\end{prob}}
\newcommand{\bede}{\begin{defin}\begin{sl}}
\newcommand{\edde}{\end{sl}\end{defin}}
\newcommand{\beco}{\begin{coro}\begin{sl}}
\newcommand{\enco}{\end{sl}\end{coro}}
\newcommand{\disp}{\displaystyle}

\newcommand{\quext}{\quad\text}
\newcommand{\qquext}{\qquad\text}
\newcommand{\de}{\partial}

\newcommand{\RR}{\mathbb{R}}

\newcommand{\EE}{\mathbb{E}}

\newcommand{\NN}{\mathbb{N}}

\newcommand{\beeq}[1]{\begin{equation}\label{#1}}
\newcommand{\eddeq}{\end{equation}}

\newcommand{\beeqa}[1]{\begin{eqnarray}\label{#1}}
\newcommand{\eddeqa}{\end{eqnarray}}

\newcommand{\beal}[1]{\begin{align}\label{#1}}
\newcommand{\eddal}{\end{align}}

\newcommand{\bespl}[1]{\begin{split}\label{#1}}
\newcommand{\edspl}{\end{split}}

\newcommand{\bega}[1]{\begin{gather}\label{#1}}
\newcommand{\edga}{\end{gather}}

\newcommand{\beeqax}{\begin{eqnarray*}}
\newcommand{\eddeqax}{\end{eqnarray*}}

\def\qed{\ifmmode % if math mode, assume display: omit penalty etc.
  \else \leavevmode\unskip\penalty9999 \hbox{}\nobreak\hfill
  \fi
  \quad\hbox{\hskip.5em\vrule width.4em height.6em depth.05em\hskip.1em}}
\def\endproofsym{\qed}
\renewenvironment{proof}[1][Proof]{\trivlist\item[\hskip\labelsep{\hskip0pt
    {\normalfont\scshape#1.}\hskip .321429\parindent}]\ignorespaces}
{\endproofsym\endtrivlist}
\def\endnobox{\def\endproofsym{}\end{proof}\def\endproofsym{\qed}}

\newcommand{\no}{\nonumber}

\newcommand{\beeqao}{\begin{eqnarray}\no}
\newcommand{\bealo}{\begin{align}\no}
\newcommand{\besplo}{\begin{split}\no}
\newcommand{\begao}{\begin{gather}\no}

\newcommand{\duav}[1]{\langle{#1}\rangle}

\newcommand{\duavb}[1]{\big\langle{#1}\big\rangle}

\newcommand{\perogni}{\forall\,}

\newcommand{\itt}{\int_0^t}

\newcommand{\io}{\int_\Omega}

\newcommand{\iTT}{\int_0^T}

\newcommand{\epsi}{\varepsilon}
\newcommand{\ee}{_{\varepsilon}}
\newcommand{\eeN}{_{\epsilon,N}}

\newcommand{\znnk}{_{0,n,k}}

\newcommand{\OO}{_{\Omega}}

\newcommand{\lhs}{left hand side}
\newcommand{\rhs}{right hand side}

\newcommand{\devv}{\partial_{V,V'}}
\newcommand{\devvs}{\partial_{V',V}}

\DeclareMathOperator{\deriv}{d}

\DeclareMathOperator{\sign}{sign}

\DeclareMathOperator{\loc}{loc}

\newcommand{\HUH}{H^1(0,T;H)}

\newcommand{\HUVp}{H^1(0,T;V')}

\newcommand{\LDH}{L^2(0,T;H)}
\newcommand{\LDV}{L^2(0,T;V)}
\newcommand{\LDVp}{L^2(0,T;V')}

\newcommand{\LIV}{L^\infty(0,T;V)}
\newcommand{\LIVp}{L^\infty(0,T;V')}

\newcommand{\LDHD}{L^2(0,T;H^2(\Omega))}

\let\TeXchi\chi
\def\chi{{\setbox0 \hbox{\mathsurround0pt
$\TeXchi$}\hbox{\raise\dp0 \copy0 }}}

\newcommand{\znn}{_{0,n}}

\newcommand{\zzn}{_{0,n}}

\newcommand{\teta}{\vartheta}

\newcommand{\bhat}{\widehat{b}}
\newcommand{\bciapo}{\widehat{b}}

\newcommand{\calA}{{\mathcal A}}

\newcommand{\calE}{{\mathcal E}}

\newcommand{\calJ}{{\mathcal J}}

\newcommand{\tetasu}{\overline{\teta}}
\newcommand{\tetagiu}{\underline{\teta}}

\newcommand{\dit}{\deriv\!t}
\newcommand{\dis}{\deriv\!s}

\newcommand{\diy}{\deriv\!y}

\newcommand{\ddt}{\frac{\deriv\!{}}{\dit}}

\numberwithin{equation}{section}
\begin{document}

\title{Asymptotic uniform 
 boundedness of energy solutions to
 the Penrose-Fife model}

\author{Giulio Schimperna\\
Dipartimento di Matematica, Universit\`a di Pavia,\\
Via Ferrata~1, I-27100 Pavia, Italy\\
E-mail: {\tt giusch04@unipv.it}\\
\and
Antonio Segatti\\
Dipartimento di Matematica, Universit\`a di Pavia,\\
Via Ferrata~1, I-27100 Pavia, Italy\\
E-mail: {\tt antonio.segatti@unipv.it}\\
\and
Sergey Zelik\\
Department of Mathematics, University of Surrey,\\
Guildford, GU2 7XH, United Kingdom\\
E-mail: {\tt S.Zelik@surrey.ac.uk}
}

%\date{}
\maketitle
\begin{abstract}
 We study a Penrose-Fife phase transition model coupled
 with homogeneous Neumann boundary conditions. Improving previous
 results, we show that the initial value problem 
 for this model admits a unique solution under 
 weak conditions on the initial data.
 Moreover, we prove asymptotic regularization properties
 of weak solutions.
\end{abstract}


\vspace{2mm}

\noindent {\bf AMS (MOS) subject clas\-si\-fi\-ca\-tion:}~~%
 35B40, 35K45, 80A22.

%%%%%%%%%%%%%%%%%%%%%%%%%%%%%%%%%%%%%%%%%%%%%%%%%%%%%%%%%%%%%%%%%%%
%%%%%%%%%%%%%%%%%%%%%% Caption da ambo i lati %%%%%%%%%%%%%%%%%%%%%

%
%
%\pagestyle{myheadings}
%\newcommand\testopari{\sc Giulio~Schimperna}
%\newcommand\testodispari{\sc Cahn-Hilliard equation}
%\markboth{\testodispari}{\testopari}

%%%%%%%%%%%%%%%%%%%%%%%%%%%%%%%%%%%%%%%%%%%%%%%%%%%%%%%%%%%%%%%%%%%
%%%%%%%%%%%%%%%%%%%%%%%%%%%%%%%%%%%%%%%%%%%%%%%%%%%%%%%%%%%%%%%%%%%

\section{Introduction}

The Penrose-Fife system, proposed by O.~Penrose and P.~Fife in \cite{PF1,PF2},
represents a thermodynamically consistent model for the description 
of the kinetics of phase transition and phase separation processes 
in binary materials. It couples the singular heat equation \eqref{calore} 
for the {\sl absolute}\/ temperature $\teta$ with 
a nonlinear relation describing the evolution of the {\sl phase variable}\/ 
$\chi$ which represents the local proportion of one of the two components.
This can be of the {\sl fourth}\/ order in space
(cf.~\eqref{phase1}-\eqref{phase2} below), 
in case the physical process preserves the total mass of $\chi$ 
({\sl conserved}\/ Penrose-Fife model, describing phase separation) 
or of the {\sl second}\/ order in space (cf.~\eqref{phase} below), 
in case the total mass of $\chi$ is admitted to 
vary ({\sl non-conserved}\/ Penrose-Fife model,
describing phase transition).
In the conserved case, the equation for $\chi$
is usually written as a system by introducing an auxiliary
variable $w$ called {\sl chemical potential}.
We refer to the next section for a detailed presentation of
the equation and to the papers \cite{CGRS,CL,SZ} for further
mathematical background.

Due to physical considerations, it is generally accepted to
consider no-flux boundary conditions for $\chi$ and,
in the conserved case, also for $w$.
On the other hand, various types of boundary conditions 
(for instance, no-flux, non-homogeneous Dirichlet, or Robin
conditions) make sense for the heat equation \eqref{calore}, 
which give rise to different mathematical scenarios.
We refer the reader to \cite{CL,CLS,DK,La1,La2,SZ} for the 
case of Robin (or ``third type'') conditions, to 
\cite{FS,GM} for the Dirichlet case, and to \cite{CGRS,IK,KK} for
the homogeneous Neumann case, which is probably the most
difficult one due to lower coercivity properties of the
elliptic operator in \eqref{calore}.

Our aim in this paper is that of improving existing results 
on the homogeneous Neumann problem for the Penrose-Fife model
both in the non-conserved and in the conserved case.
Actually, our results will cover both situations,
generally with minor variations in the proofs.
As a first property, we will show that the problem admits 
a unique solution under weak assumptions on the initial data.
Actually, noting that the system admits 
a natural Liapounov functional representing the total energy,
we will prove that a (unique) weak solution exists for any initial data 
$(\teta_0,\chi_0)$ having finite energy.
 This improves existing results which assume some extra summability
condition on $\teta_0$, typically $\teta_0\in L^2(\Omega)$. 
As a side effect, we pay the price that the heat equation has
to be interpreted in the generalized $(H^1)'$-framework 
developed by Damlamian and Kenmochi in
\cite{DK}. Namely, the (thermal part of the) energy has
to be intended as a (relaxed) functional operating on
the negative order Sobolev space $H^1(\Omega)'$
(cf.~\eqref{defiE} below)
and also the relation linking $\teta$ to its inverse
has to be stated properly.
On the other hand, if we know in addition that
$\teta_0\in L^1(\Omega)$, then we can prove that
$\teta(t) \in L^1(\Omega)$ for all $t\ge 0$; 
moreover, both the energy functional and the 
relation between $\teta$ and its inverse can 
be written in the usual (pointwise) sense
(cf.~\eqref{defiEstrong} below). 

In the subsequent part of the paper, we prove our main results,
which regard uniform time-regularization properties of weak solutions.
In this frame, we will actually present two theorems. Firstly, we will
show that, for any $T>0$, there exists a constant $\tetagiu>0$ depending
only on the ``energy'' of the initial data and on $T$ such that 
$\teta(x,t)\ge \tetagiu$ for a.e.~$(x,t)\in(T,+\infty)\times\Omega$.
Second, we will prove that a similar bound from above
(i.e., $\teta(x,t)\le \tetasu$ for a.e.~$(x,t)\in(T,+\infty)\times\Omega$)
and a suitable $\tetasu>0$) holds provided that the initial temperature
$\teta_0$, in addition to the ``energy'' regularity, satisfies 
the additional hypothesis $\teta_0\in L^{3+\epsi}(\Omega)$ for some
$\epsi>0$.  This additional condition appears also in other works concerning
$L^\infty$-regularization properties of the solutions to {\sl very fast diffusion}\
equations like \eqref{calore} on $\mathbb{R}^3$ 
(see, e.g., \cite{vazlibro}, \cite{bV-Adv} and the references therein).
In particular, in three space dimensions, the exponent $p=3$ happens to be critical 
for the boundedness (for strictly positive times) of the solutions to \eqref{calore}: 
starting from initial data in $L^p(\mathbb{R}^3)$, $p>3$, implies boundedness of the 
solutions for strictly positive times (see \cite{bV-Adv}). For $p<3$, 
the situation is drastically different, as the self similar solution \eqref{similarity} shows.

We finally note that the proved uniform bounds permit, by standard
methods, to improve furtherly the regularity of solutions for strictly
positive times. In particular, our estimates complement a recent 
paper by Pr\"uss and Wilke \cite{PW} who show maximal regularity estimates
for the conserved Penrose-Fife model under the conditional (i.e., 
unproved in their paper) assumption that the temperature 
$\teta$ satisfies the uniform bounds $\tetagiu \le \teta(x,t)\le \tetasu$ 
almost everywhere. Thanks to our results, the maximal regularity estimates
of Pr\"uss and Wilke hold for all strictly positive times and all 
weak solutions emanating from initial data satisfying the
``energy regularity'' plus the condition $\teta_0\in L^{3+\epsi}(\Omega)$.
It is also worth noting that the uniform bound $\teta\ge \tetagiu$ 
implies, in the nonconserved case, the so-called ``separation 
from singularities'' property for $\chi$
in case the configuration potential of the system 
(i.e., the function $\bhat$ defined in \eqref{defibhat})
has a bounded domain. For instance, in the physically relevant
case of the {\sl logarithmic potential}
\begin{equation}\label{potlog}
  \bhat(r)=(1+r)\log(1+r)+(1-r)\log(1-r),
\end{equation}
whose domain is $[-1,1]$, this means that 
$-1+\delta\le\chi(x,t)\le1-\delta$ for some $\delta>0$ 
and for all times $t\ge 1$. Unfortunately, we cannot prove
such a property in the case of the conserved model. Actually,
up to our knowledge, this is unknown, at least in
the three-dimensional case, also in for the (simpler)
conserved phase-field model of Caginalp type (cf.~the 
related discussion in~\cite{MZ}).

The remainder of the paper is organized as follows. In the next section
we will present our hypotheses and state our results. The proofs
will be detailed in the subsequent Section~\ref{sec:proofs}.

\smallskip

\noindent%
{\bf Acknowledgment.}~~The authors would like to thank Prof.~Philippe
Lauren\c cot for discussions regarding the strategy of some proofs.

%%%%%%%%%%%%%%%%%%%%%%%%%%%%%%%%%%%%%%%%%%%%%%%%%%%%%%%%%%%%%%%%%%%
%%%%%%%%%%%%%%%%%%%%%%%%%%%%%%%%%%%%%%%%%%%%%%%%%%%%%%%%%%%%%%%%%%%

\section{Main results}
\label{sec:main}

Let $\Omega$ be a smooth bounded domain of $\RR^d$, 
$d\in\{2,3\}$. For the sake of simplicity,
let us assume $|\Omega|=1$ so that 
$\| v \|_{L^p(\Omega)}\le \| v \|_{L^q(\Omega)}$
for all $1\le p \le q \le +\infty$, $v\in L^q(\Omega)$.
For simplicity, we will often write 
$\| \cdot \|_p$ in place of $\| \cdot \|_{L^p(\Omega)}$.
Let $H:=L^2(\Omega)$, endowed with the standard scalar 
product $(\cdot,\cdot)$ 
and norm $\| \cdot \|$. Let also $V:=H^1(\Omega)$. 
We note by $\|\cdot\|_X$
the norm in the generic Banach space $X$.
and by $\langle\cdot,\cdot\rangle_X$
the duality between $X'$ and $X$.

For any function, or functional $z$, defined on
$\Omega$, we can then set
\begin{equation}\label{defim}
  z\OO:=\frac1{|\Omega|}\io z
   = \io z,  
\end{equation}
where the integral is substituted with 
the duality $\langle z,1 \rangle$ in case, e.g.,
$z\in V'$. 

We also define the elliptic operator
\begin{equation}\label{defiA}
  A:V\to V',\qquad
   \duav{Av,z}:=\io \nabla v\cdot \nabla z.
\end{equation}
Then, for $t\in(0,+\infty)$, we consider the 
singular heat equation
\begin{align}\label{calore}
  & \teta_t + A u = - \chi_t,
    \quext{in }\,V',\\
 \label{identif}
  & u = -\frac1\teta,
   \quext{almost everywhere in }\,\Omega.
\end{align}
In the {\sl non-conserved}\/ case, this is
coupled with the equation
\begin{equation}\label{phase}
  \chi_t + A \chi + b(\chi) - \chi = u,
  \quext{in }\,H,
\end{equation}
while in the {\sl conserved}\/ case, 
\eqref{calore}-\eqref{identif}
is coupled with the system
\begin{align}\label{phase1}
  & \chi_t + A w = 0,
  \quext{in }\,V',\\
 \label{phase2}
  & w = A \chi + b(\chi) - \chi - u,
  \quext{in }\,V'.
\end{align}
In both cases, we will take
the initial conditions
\begin{equation} \label{iniz}
  \teta|_{t=0}=\teta_0, \quad \chi|_{t=0}=\chi_0, 
   \quext{in }\,\Omega.
\end{equation}
The nonlinear function $b$ is assumed to satisfy
\begin{equation}\label{hpb}
  b\in C^{0,1}_{\loc}(I;\RR), \quad b(0)=0,
   \quad b'(r)\ge 0~~\text{a.e.~in }\,I,
   \quad \lim_{r\to\de I} b(r)\sign r=\lim_{r\to\de I} b'(r)=+\infty,
\end{equation}
where $I$, the {\sl domain}\/ of $b$, is an open, possibly
bounded, interval of $\RR$ containing $0$.
We also set
\begin{equation}\label{defibhat}
  \bhat(r):=\int_0^r b(s)\,\dis
\end{equation}
in such a way that $\bhat$ is a convex function 
satisfying $\bhat(0)=0$. 

Concerning the initial datum $\chi_0$, we will
always assume
\begin{equation}\label{hpchi0}
  \chi_0\in V, \qquad \bhat(\chi_0)
   \in L^1(\Omega).
\end{equation}
To specify the regularity required for the initial 
temperature $\teta_0$, we first need to introduce some
convex analysis machinery. First of all, we set
\begin{equation}\label{defij}
  j(v):= - \log (- v),
   \qquad j^*(z):= - 1 - \log z,
\end{equation}
%
%{\bf Controlla. Se \`e giusto, la costante 1 dentro $j^*$ va ovunque 
%cambiata in $-1$.}
%
respectively for $v<0$ and $z>0$. Then, $j$ and
$j^*$ can be seen as a couple of {\sl conjugate
functions}\/ according to the standard theory
(cf., e.g., ~\cite{Br}). This permits to introduce
the convex functional
\begin{equation}\label{defiJ}
  J:V\to [0,+\infty],
   \quad J(v):= \io j(v), 
\end{equation}
where we have implicitly set $J(v)=+\infty$ for those
$v\in V$ such that $j(v)$ is not summable
(this happens, for instance, when $v$ is strictly 
positive on a set of strictly positive measure).
The set $\{v\in V:~J(v)<+\infty\}$ is called 
{\sl domain}\/ of $J$. Then, the conjugate of
$J$ is given by
\begin{equation}\label{defiJstar}
  J^*:V' \to (-\infty,+\infty], \qquad
   J^*(\zeta):= \sup \big\{ \duav{\zeta,v} - J(v),~v\in V\big\}.
\end{equation}
We can also introduce the subdifferential of $J$ with 
respect to the duality between $V'$ and $V$. 
Namely, given $v\in V$, we set
\begin{equation}\label{weaksubd}
  \zeta\in \devv J(v)\Longleftrightarrow
   \duav{\zeta,z-v}\le J(z)-J(v)
   \quad \perogni z\in V.
\end{equation}
In general, the ``weak'' subdifferential $\devv J$
is a multivalued maximal-monotone operator. Its 
structure is analyzed in several
papers (see, e.g., \cite{BCGG}; see also \cite{brezisart}
and \cite{GrR} for the slightly different situation
where $V$ is substituted by $H^1_0(\Omega)$).
Actually, it in not difficult to prove that
\begin{equation}\label{weaksubd2}
  \zeta\in L^1(\Omega),~~ \zeta = j'(v) = - \frac1v~~
   \text{a.e.~in }\,\Omega \Longrightarrow
   \zeta\in \devv J(v).
\end{equation}
On the other hand, for given $v\in V$, the set 
$\devv J(v)$ needs not be contained in $L^1(\Omega)$.
More precisely, a generic element $\zeta\in\devv J(v)$ 
is a bounded measure which may have a singular part
$\zeta_s$. A precise characterization is given in 
\cite[Thm.~3]{brezisart} for the $H^1_0$-case. However,
it is easy to realize that such a characterization
extends to the case of $V=H^1(\Omega)$ at least provided
$\Omega$ is smooth, the only difference being that 
$\zeta_s$ can be supported also on the boundary $\de\Omega$.
On the other hand, if we are able to prove that
the singular part of some element $\zeta$ is $0$,
then we still have pointwise inclusion. Namely,
there holds that
\begin{equation}\label{weaksubd3}
  \zeta \in \devv J(v)\cap L^1(\Omega)
   \Longrightarrow \zeta = - \frac1v~~
   \text{a.e.~in }\,\Omega.
\end{equation}
That said, the minimal regularity required
on $\teta_0$ is given by
\begin{equation}\label{hpteta0}
  \teta_0\in V', \quad
   J^*(\teta_0) < +\infty.
\end{equation}
Then, we can define the ``strong'' energy functional
of the system as
\begin{equation}\label{defiEstrong}
  \calE_s(\teta,\chi):=
   \io\Big(\teta + \chi - 1 - \log\teta
     + \frac12|\nabla\chi|^2 + \bhat(\chi)- \frac12\chi^2 \Big),
\end{equation}
However, since in general $\teta_0$ is just an element
of $V'$, it does not make sense to compute $\calE_s(\teta_0,\chi_0)$.
Actually, we have to relax $\calE_s$, by defining
\begin{equation}\label{defiE}
  \calE(\teta,\chi):= \duav{\teta,1} + J^*(\teta)
   + \io\Big(\chi + \frac12|\nabla\chi|^2 + \bhat(\chi)- \frac12\chi^2 \Big).
\end{equation}
Actually, assumptions \eqref{hpteta0} and \eqref{hpchi0} 
are equivalent to asking that the ``relaxed'' energy 
$\EE_0:=\calE(0)$ is finite.
In the sequel we shall often write
$\calE(t)$ in place of $\calE(\teta(t),\chi(t))$.
Notice that the last \eqref{hpb} ensures the 
coercivity of $\calE$. 
Notice also that the no-flux conditions 
entail some conservation properties.
Actually, for the conserved system, testing \eqref{calore}
and \eqref{phase1} by $1$, we immediately get
\begin{equation}\label{consmedie}
  \teta\OO(t)=(\teta_0)\OO, \quad
   \chi\OO(t)=(\chi_0)\OO, \quad
   \perogni t\ge 0,
\end{equation}
while of course for the non-conserved system we only
have
\begin{equation}\label{consmediew}
  \teta\OO(t)+\chi\OO(t)=(\teta_0)\OO+(\chi_0)\OO, \quad
   \perogni t\ge 0.
\end{equation}
The following result, stating existence of
strong solutions, is well-known
(see, for instance, \cite{KK} for a proof):
\bete\label{teo:esi:pf}
 Let us assume\/ \eqref{hpb}, 
 \eqref{hpteta0}, \eqref{hpchi0} and,
 additionally, let
 \begin{equation}\label{tetazl2}
   \teta_0 \in H.
 \end{equation}
 Then, in the\/ {\rm non-conserved} case,
 there exists a unique triplet $(\teta,u,\chi)$
 satisfying, for all $T>0$,
 \begin{align}\label{regoteta}
   & \teta-\log\teta \in L^\infty(0,T;L^1(\Omega)), \quad
    \teta\in L^2(0,T;H),\\
  \label{regou}
   & u \in L^2(0,T;V),\\ 
  \label{regochi}
   & \chi \in \HUH \cap \LIV \cap\LDHD,
    \quad b(\chi) \in \LDH,
 \end{align}
 and solving\/ \eqref{calore}-\eqref{identif} and 
 \eqref{phase} a.e.~in~$(0,T)$, together with the initial 
 conditions\/ \eqref{iniz}. 
 Analogously, in the\/ {\rm non-conserved} case,
 we have a unique triplet $(\teta,u,\chi)$
 satisfying, for all $T>0$, \eqref{regoteta}-\eqref{regou},
 together with
 \begin{equation}\label{regochiw}
   \chi \in \HUVp \cap \LIV \cap\LDHD,
    \quad b(\chi) \in \LDH,
 \end{equation}
 solving\/ \eqref{calore}-\eqref{identif} and 
 \eqref{phase1}-\eqref{phase2} 
 a.e.~in~$(0,T)$, and satisfying\/ \eqref{iniz}. 
\ente
\beos\label{rem:Jstarmeglio}
 It is worth noting (cf.~\cite{BCGG} for more details) that 
 \eqref{tetazl2} entails 
 %
% \begin{equation}\label{teta0meglio}
%   J^*(\teta_0) = \io ( 1 - \log \teta_0 ).
% \end{equation}
 % 
%{\asnote piuttosto, non dovrebbe essere 
\begin{equation}\label{teta0meglio}
   J^*(\teta_0) = \io ( -1 - \log \teta_0 ).
 \end{equation}
%grazie a \eqref{defij} e \eqref{tetazl2}??
%}
 Thus, in fact no functional on $V'$ appears in the 
 above ``strong'' formulation. 
\eddos
\noindent%
In the next result we will provide an existence 
theorem working without the additional regularity 
\eqref{tetazl2}. We will pay the price of the occurrence
of relaxed functionals. Moreover, we will also need to
intend relation \eqref{identif} in the relaxed sense
of \cite{DK}.
\begin{teor}[$\boldsymbol{(V')}$-solutions]
 \begin{sl}\label{teo:esi:w}
 Let\/ \eqref{hpb}, \eqref{hpteta0} and \eqref{hpchi0} hold.
 %and let also
 %
 %\begin{equation}\label{tetazVp}
 %  \teta_0 \in V'.
 %\end{equation}
 %
 Then, there exists a unique triplet $(\teta,u,\chi)$, satisfying,
 for all $T>0$,
 \begin{equation}\label{regotetaw}
   \teta\in \HUVp, \quad
    J^*(\teta) \in L^\infty(0,T),
 \end{equation}
 together with\/ \eqref{regou} and, a.e.~in~$(0,T)$, equation 
 \eqref{calore}. Moreover, in the non-conserved case,
 \eqref{regochi} and, a.e.~in~$(0,T)$, equation \eqref{phase} 
 hold, while, in the conserved case, \eqref{regochiw}
 and, a.e.~in~$(0,T)$, equations \eqref{phase1}-\eqref{phase2}
 hold. Moreover, there hold the initial condition\/ \eqref{iniz} and,
 a.e.~in~$(0,T)$, the weak identification property
 \begin{equation}\label{identifw}
   \teta \in \devv J(u).
 \end{equation} 
 Finally, any $V'$-solutions satisfies, for all $t\in [0,T]$, 
 the\/ {\rm energy equality}
 \begin{align}\no
   & \itt \big( \| \nabla u \|^2
    + \| \chi_{t} \|^2 \big)
    = J^*(\teta_0) - J^*(\teta(t)) 
    + \duav{\teta_0,1} - \duav{\teta(t),1}\\
  \label{energy-teo}
    & \mbox{}~~~~~
    + \io\Big(\chi_0 + \frac12|\nabla\chi_0|^2 + \bhat(\chi_0)- \frac12\chi_0^2 \Big)
    - \io\Big(\chi(t) + \frac12|\nabla\chi(t)|^2 + \bhat(\chi(t))- \frac12\chi(t)^2 \Big),
 \end{align}
 in the non-conserved case. In the conserved case, the same
 holds provided that the term $\|\chi_t\|^2$ on the \lhs\
 is replaced by $\| \nabla w \|^2$. 
\ente
\noindent%
We now turn to discussing regularization properties
of solutions. The first result regards
the function $u$.
\bete\label{teo:asy:u}
 Let\/ \eqref{hpb}, \eqref{hpteta0}
 and \eqref{hpchi0} hold. Then, the $V'$-solution
 either to the non-conserved or to the conserved
 problem satisfies
 \begin{align}\label{asym-u}
   & \| u \|_{L^\infty(1,\infty;V)} + \| u \|_{L^\infty((2,\infty)\times \Omega)}
    \le Q( \EE_0 ),\\
  \label{asym-chi}
   & \| \chi \|_{L^\infty(1,\infty;H^2(\Omega))} + \| b(\chi) \|_{L^\infty(1,\infty;H)}
    \le Q( \EE_0 ).
 \end{align}
 Moreover, in the\/ {\rm non-conserved} case we also have 
 the ``separation property''
 \begin{equation}\label{separ}
  \| b(\chi) \|_{L^\infty((3,\infty)\times \Omega)}
   \le Q( \EE_0 ).
 \end{equation}
 Here and below, $Q$ is a computable nonnegative function, monotone
 in each of its arguments, whose expression is independent 
 of initial data and of time.
\ente
\beos
 It is worth noting that the above regularization properties are
 in fact {\sl instantaneous}. Indeed, with minor
 modification in the proofs one could easily see that 
 \eqref{asym-u}-\eqref{asym-chi} and \eqref{separ}
 hold starting from any $\tau>0$ (and not only from $\tau=1$
 or 2 or 3).
 Of course, then the functions $Q$ on the \rhs s would also
 monotonically depend on $\tau^{-1}$ and possibly explode
 for $\tau\searrow 0$. The same considerations hold also for 
 what is proved in Theorem~\ref{teo:asy:teta} below.
\eddos
\noindent%
In case the initial temperature, beyond satisfying
\eqref{hpteta0}, is an $L^1$-function, we can say
something more precise on regularity:
\begin{teor}[$\boldsymbol{(V'\cap L^1)}$-solutions]
 \begin{sl}\label{teo:esi:w2}
 Let\/ \eqref{hpb}, \eqref{hpteta0} and \eqref{hpchi0} hold. 
 Let also 
 \begin{equation}\label{tetazLu}
   \teta_0 \in L^1(\Omega).
 \end{equation}
 Then, the triplet $(\teta,u,\chi)$
 given by\/ {\rm Theorem~\ref{teo:esi:w}}
 additionally satisfies, for all $T>0$,
 \begin{equation}\label{tetaLu}
   \teta\in C^0([0,T];L^1(\Omega))
 \end{equation}
 and the strong identification 
 property~\eqref{identif}.
\ente
\noindent%
It is easy to show that, if $\teta_0\in L^p(\Omega)$ for 
$p>1$, then $\teta(t)$ remains in $L^p(\Omega)$ for $t>0$.
What is more interesting is that, if $p>3$, then 
$\teta$ is asymptotically uniformly bounded:
\bete\label{teo:asy:teta}
 Let\/ \eqref{hpb}, \eqref{hpteta0} and
 \eqref{hpchi0} hold and let also
 \begin{equation}\label{tetazl3}
   \teta_0 \in L^{3+\epsi}(\Omega)
   \quext{for some }\,\epsi>0.
 \end{equation}
 Then, the $(V'\cap L^1)$-solution 
 either to the non-conserved or to the conserved
 problem satisfies the additional bound
 \begin{equation}\label{asym-teta}
   \| \teta \|_{L^\infty((2,\infty)\times \Omega)}
    \le Q\big( \EE_0, \| \teta_0 \|_{3+\epsi} ).
 \end{equation}
\ente
\beos
 Relation \eqref{tetazl3} suggests that $3$ should
 play the role of a critical exponent for equation
 \eqref{calore} in space dimension $3$. Actually, it is easy to 
 check that the related ``very-fast diffusion'' equation
 \begin{equation}\label{veryfast}
   \teta_t + \Delta \teta^{-1} = 0
 \end{equation}
 over the whole space $(0,+\infty)\times \RR^3$ admits 
 the similarity solution (see, e.g., \cite{vazlibro})
 \begin{equation}\label{similarity}
   \teta(t,x) = \frac{2(T-t)_+^{1/2}}{|x|},
 \end{equation}
 which belongs to $L^p_{\loc}(\RR^3)$ for all $p<3$ and
 all $t\ge 0$ and does not exhibit any instantaneous 
 regularizing effect (of course, it satisfies a delayed
 regularization property since it extinguishes in a finite time; 
 however, this effect is not expected to hold in the case 
 of a finite domain when we have conservation of mass).
 However, we do not know what happens in the critical case of
 an initial datum $\teta_0$ belonging to $L^3(\Omega)$. 
\eddos

%%%%%%%%%%%%%%%%%%%%%%%%%%%%%%%%%%%%%%%%%%%%%%%%%%%%%%%%%%%%%%%%%%%
%%%%%%%%%%%%%%%%%%%%%%%%%%%%%%%%%%%%%%%%%%%%%%%%%%%%%%%%%%%%%%%%%%%

\section{Proofs}
\label{sec:proofs}

All proofs will be in principle given only for the {\sl conserved}\/
case which is, actually, more difficult. The properties holding only
for the non-conserved case (as well any significant differences in 
the proofs) will be remarked on occurrence.

%%%%%%%%%%%%%%%%%%%%%%%%%%%%%%%%%%%%%%%%%%%%%%%%%%%%%%%%%%%%%%%%%%%

\subsection{Proof of Theorem~\ref{teo:esi:w}}
\label{sec:vp}

We start by proving existence for weak initial data.
Given $\teta_0$ satisfying \eqref{hpteta0},
we then set, for $n\in\NN$,
\begin{equation}\label{tetaznn}
  \teta\zzn + n^{-1} A \teta\zzn = \teta_0.
\end{equation}
The properties of this approximation deserve to be 
stated in a lemma.
\bele\label{lemma:znn}
 Let \eqref{hpteta0} hold and let 
 $\teta\znn$ be defined by\/ \eqref{tetaznn}. 
 Then, $\teta\znn\in V$ for all $n\in\NN$. Moreover,
 \begin{equation}\label{tetaznn2}
   \teta\zzn\to \teta_0\quext{strongly in }\,V'\,\,\,\,\hbox{ and } \,\,\,\,
    J^*(\teta\zzn) \le J^*(\teta_0)\,\,\,\forall\,n\in \mathbb{N}. 
 \end{equation}
 Moreover, if also\/ \eqref{tetazLu} holds, then we also
 have
 \begin{equation}\label{tetaznn2b}
   \teta\zzn\to \teta_0\quext{strongly in }\,L^1(\Omega)
 \end{equation}
 and
 \begin{equation}\label{tetaznn3}
   \log \teta\znn \in L^1(\Omega), \qquad
    \io - 1 - \log \teta\znn 
    \le \io -1 - \log \teta_0,
    \quad\perogni n\in\NN.
 \end{equation}
\enle
\begin{proof}
Being $\teta_0\in V'$, it is clear that $\teta\znn\in V$ 
for all $n$. Moreover, the $V'$-strong convergence in 
\eqref{tetaznn2} can be proved by standard Hilbert techniques.
In order to complete the proof of \eqref{tetaznn2}, we introduce,
for any given $n\in \mathbb{N}$, the strictly positive sequence 
$\teta\znnk:=\max\{\teta\znn,\frac{1}{k}\}$, for $k\in \mathbb{N}$. 
Correspondingly, we set $u\znnk := -\frac{1}{\teta\znnk}$. Note that, 
by construction, $u\znnk \in V$ and the map
$\teta\znn\mapsto\teta\znnk$ is monotone. 
As a consequence, we can write
\begin{equation}\label{antonio11}
  0 \ge \Big\langle -\frac{1}{n}A \teta\znn, u\znnk \Big\rangle
   = \langle \teta\znn-\teta_0,u\znnk \rangle
   = \langle \teta\znnk-\teta_0,u\znnk \rangle 
       + \langle \teta\znn-\teta\znnk, u\znnk \rangle.
\end{equation}
Now, since $\teta\znnk = -\frac{1}{u\znnk}$ almost everywhere 
in $\Omega$, we have $\teta\znnk \in \partial_{V,V'} J(u\znnk)$.
Equivalently, $u\znnk \in \partial_{V',V} J^*(\teta\znnk)$ 
(where the subdifferential acts now in the duality between $V'$ and $V$)
for any $k$. Hence,
\begin{equation}\label{antonio12}
  \langle \teta\znnk-\teta_0,u\znnk \rangle
   \ge J^*(\teta\znnk) - J^*(\teta_0),
\end{equation}
by definition of subdifferential. 
On the other hand, since $\teta\znnk \equiv \teta\znn$ in 
$\Omega\cap \{ \teta\znn\ge 1/k  \}$,
\begin{equation}\label{antonio13}
  \langle \teta\znn-\teta\znnk,u\znnk \rangle
   = \int_{\Omega\cap \{\teta\znn \le 1/k \}}
   (-k) \cdot (\teta\znn - \teta\znnk) \ge 0.
\end{equation}
Thus, collecting the above computations we have
\begin{equation}\label{antonio14}
 J^*(\teta\znnk)\,\le\, J^*(\teta_0).
\end{equation}
Finally, since for any $n\in\mathbb{N}$ we have that
$\teta\znnk\xrightarrow{k\nearrow +\infty}\teta\znn$
strongly in $L^p(\Omega)$ for any $p\in[1,6)$ 
(hence, a fortiori, in $V'$), 
we have the following chain of inequalities:
\begin{equation}\label{antonio15}
  J^*(\teta\znn) 
   \le \liminf_{k\nearrow +\infty} J^*(\teta\znnk)
   \le \limsup _{k\nearrow +\infty} J^*(\teta\znnk)
   \le  J^*(\teta_0), \quext{ for any }\, n\in\mathbb{N},
\end{equation}
i.e., \eqref{tetaznn2} holds.

Now, we assume that also \eqref{tetazLu} holds, 
namely we assume that $\teta_0\,\in\, V´\cap L^1(\Omega)$. 
Note that, thanks to Remark~\ref{rem:Jstarmeglio},
in this new regularity framework, \eqref{tetaznn3} 
is nothing else than $J^*(\teta\znn) \le J^*(\teta_0)$,
which as been proved above.

Thus, we only need to prove the $L^1$-convergence. 
To this end, we have to be a bit more careful. 
First, we define the Banach space $X:=V'\cap L^1(\Omega)$,
endowed with the norm 
$$
  \| \cdot \|_X
   := \| \cdot \|_{V'} + \| \cdot \|_{L^1(\Omega)},
$$
and introduce the unbounded linear operator $\calA$ on $X$ 
defined as $\calA v:= Av$ with domain 
$$
  D(\calA):=\big\{ v \in V:~\Delta v \in L^1(\Omega) \big\},
$$
where $\Delta$ is the usual distributional Laplace operator.
Then, we have that $\calA$ is an accretive operator
on the space $X$. Indeed, by \cite[Prop.~II.3.1]{Ba}, 
this corresponds to checking 
that, if $\lambda>0$ and
\begin{equation}\label{accre11}
  x_i + \lambda \calA x_i = f_i, \quad i=1,2,
\end{equation}
for $f_i\in X$, then 
\begin{equation}\label{accre12}
  \| x_1 - x_2 \|_X \le \| f_1 - f_2 \|_X.
\end{equation}
Actually, the analogue of \eqref{accre12} w.r.t.~the $V'$-norm
can be obtained by testing the difference
\begin{equation}\label{accre13}
  x_1 - x_2 + \lambda \calA ( x_1 - x_2) = f_1 - f_2
\end{equation}
by $(I+A)^{-1}(x_1 - x_2)$ where $I$ is the identity mapping of $H$
(and, hence, $I+A:V\to V'$ is the Riesz isomorphism).
On the other hand, the $L^1$-analogue of \eqref{accre12}
is obtained by testing \eqref{accre13} by $\sign(x_1-x_2)$ and
applying the Brezis-Strauss theorem \cite[Lemma~2]{BS}.
Moreover, we have that $D(\calA)$ is dense in $V'\cap L^1(\Omega)$.
To see this, let us take $z\in V'\cap L^1(\Omega)$. 
Then, setting 
\begin{equation}\label{accre13b}
  z_n:= \rho_n * z = \io \rho_n(x-y) z(y) \,\diy,
\end{equation}
where $\rho_n$ is the standard mollifer,
it is clear that $z_n$ is smooth (hence, in particular,
it belongs to $D(A)$). Moreover, the convergence 
$z_n\to z$ in $L^1(\Omega)$ follows from standard
properties of convolutions, while the convergence
$z_n\to z$ in $V'$ follows from the density of
$H$ in $V'$ and from the fact that the mapping
$z\mapsto z_n$ is a contraction w.r.t.~the $V'$-norm
(this may be verified for $z\in H$ by using Fubini's
theorem and then extended to $V'$ by density). 
These facts permit to apply \cite[Prop.~3.2 (e)]{Ba},
which gives exactly the convergence property 
\eqref{tetaznn2b}, which concludes the proof.
\end{proof}
\noindent%
Thus, taking $\teta\znn$ as an initial datum for 
equation~\eqref{calore} (while the initial datum
$\chi_0$ is kept fixed), existence of a corresponding 
solution $(\teta,u,\chi)$ % $(\teta_n,u_n,\chi_n)$ 
is guaranteed by Theorem~\ref{teo:esi:pf}. 
Our aim will be now that of 
removing the approximation of the initial datum letting 
$n\nearrow+\infty$. With this aim, we start by recalling 
a couple of basic a-priori estimates.
The procedure is detailed only in the conserved case,
the differences occurring in the non-conserved case being
pointed out at the end. For the meanwhile, we will not
emphasize the dependence on $n$ in the notation.

\smallskip

\noindent%
{\bf Energy estimate.}~~%
We test \eqref{calore} by $1+u$, \eqref{phase1}
by $w$ and \eqref{phase2} by $\chi_t$.
This formal procedure will be justified at the end,
when we prove \eqref{energy-teo}. We obtain
\begin{equation}\label{energy}
  \ddt \calE + \| \nabla u \|^2
   + \| \nabla w \|^2 
   = 0,
\end{equation}
where $\calE$ was defined in \eqref{defiE}.
Using also the properties of $A$, we immediately get
\begin{equation}\label{st-energy}
  \| \calE \|_{L^\infty(0,T)}
   + \| \nabla u \|_{L^2(0,T;H)} 
   + \| \nabla w \|_{L^2(0,T;H)} 
   + \| \chi_t \|_{L^2(0,T;V')} 
  \le c \EE_0.
\end{equation}
Here and below, the letters $c$ and $\kappa$ 
will denote generic positive constants, 
independent of initial data and of time, whose 
value possibly varies on occurrence,
$\kappa$ being used in estimates from below.
In particular, the above estimate is uniform with 
respect to $T$.

\smallskip

\noindent%
{\bf A generalized Poincar\'e inequality.}~~%
To estimate the full $V$-norm of $u$ (and not just the
$H$-norm of its gradient), we need a proper form
of Poincar\'e's inequality (cf., e.g., \cite[Lemma~5.1]{Ke} 
for a similar tool), which we prove just for the sake 
of completeness:
\bele\label{lemma-log-poincare}
 Assume $\Omega$ is a bounded open subset of $\mathbb{R}^d$. Suppose
 $v\in W^{1,1}(\Omega)$ and $v\ge 0$ a.e. in $\Omega$.
 Then, setting $K:=\int_{\Omega}(\log v)^{+}$, 
 the following estimate holds:
 \begin{equation}\label{log-poincare}
   \|v\|_{L^1(\Omega)} \le |\Omega| e ^{C_1 K} 
    + \frac{C_2}{|\Omega|}\|\nabla v\|_{L^1(\Omega)},
 \end{equation}
 the constants $C_1$ and $C_2$ depending only on $\Omega$. 
\enle
\begin{proof}
 First of all, we recall that for any function $z\in W^{1,1}(\Omega)$
 such that $|E_0|>0$ (with $E_0:=\left\{x\in \Omega : z(x)=0\right\}$)
 the following Poincar\'e type inequality
 (see \cite[Lemma 5.1, pag. 89]{LSU}) holds:
 \begin{equation}\label{start}
   \|z\|_{L^{1}(\Omega)} \le \frac{C}{|E_0|}\|\nabla z\|_{L^1(\Omega)},
 \end{equation}
 where the constant $C$ can be explicitely computed 
 and depends only on $\Omega$. 
 Now, let $v$ be a function in the hypothesis of 
 the Lemma. Set $K:=\int_{\Omega}(\log v)^{+}$ 
 and note that, thanks to the Chebychev inequality, 
 we have, for any fixed $N>0$,
 \begin{equation}\label{cheb1}
   \big|\left\{x\in \Omega: (\log v)^+ >N\right\}\big| \le \frac{K}{N},
 \end{equation}
 and consequently
 \begin{equation}\label{cheb2}
   \big|\left\{ x\in \Omega: v > N \right\}\big|\le \frac{K}{\log N}.
 \end{equation}
 Thanks to \eqref{cheb1} and \eqref{cheb2}, we can 
 fix $\bar{N} = \bar{N}(K,\Omega) = e^{\frac{2K}{|\Omega|}}$
 in such a way that 
 \begin{equation}\label{cons-cheb} 
   \big|\left\{x\in\Omega : v \le \bar{N} \right\}\big| 
    = |\Omega| - \big|\left\{x\in \Omega: v > \bar{N}\right\}\big|
   \ge \frac{|\Omega|}{2}.
 \end{equation}
 As a consequence, the inequality \eqref{start}, with 
 $z = (v - \bar{N})^{+}$ and \eqref{cons-cheb},
 entails
 \begin{align}\no
   \disp \|v\|_{L^1(\Omega)} 
    & \le \bar{N}|\Omega|
     + \| (v - \bar{N})^+\|_{L^1(\Omega)}
    \le \bar{N}|\Omega|
     + \frac{C}{\big|\left\{ v\le \bar{N}\right\}\big|}
        \|\nabla( v - \bar{N})^+\|_{L^1(\Omega)}\\
  \label{final}
   & \disp \le |\Omega|e^{\frac{2K}{|\Omega|}} 
     + 2C|\Omega| \| \nabla v\|_{L^1(\Omega)},
 \end{align}
 which is \eqref{log-poincare} with $C_1= 2/|\Omega|$ 
 and $C_2 = 2C$, $C$ being the constant
 in $\eqref{start}$.
\end{proof}

\smallskip

\noindent%
{\bf Consequences of the energy estimate.}~~%
Using the above lemma, \eqref{st-energy} 
additionally gives
\begin{equation}\label{st-21}
  \| u \|_{L^2(t,t+1;V)} 
   \le Q(\EE_0) \quad\perogni t\ge 0.
\end{equation}
Applying standard techniques
to system \eqref{phase1}-\eqref{phase2}, we also have
\begin{equation}\label{st-22}
  \| \chi \|_{L^2(t,t+1;H^2(\Omega))} 
   \le Q(\EE_0).
\end{equation}
Moreover, testing \eqref{phase2} by $\chi-\chi\OO$
and proceeding, e.g., as in the Appendix of 
\cite{MZ}, it is not difficult to arrive at
\begin{equation}\label{st-23}
  \| b(\chi) \|_{L^2(t,t+1;H)} 
   \le Q(\EE_0),
\end{equation}
whence a comparison of terms in \eqref{phase2} and 
estimate \eqref{st-energy} also give
\begin{equation}\label{st-24}
  \| w \|_{L^2(t,t+1;V)} 
   \le Q(\EE_0).
\end{equation}
\smallskip

\noindent%
{\bf Passage to the limit.}~~%
We will now let $n\nearrow+\infty$, still referring to the
conserved case. With this aim, we 
rename as $(\teta_n,u_n,\chi_n)$ the solution to the $n$-approximation.
By \eqref{st-energy} and \eqref{st-21}-\eqref{st-24}, we then
have, for any $T>0$,
\begin{align}\label{coa1}
  & u_n \to u \quext{weakly in }\,\LDV,\\
 \label{coa2}
  & \chi_n \to \chi \quext{weakly in }\,\HUVp\cap\LIV\cap\LDHD,\\
 \label{coa3}
  & w_n \to w \quext{weakly in }\, \LDV,\\
 \label{coa4}
  & b(\chi_n) \to b \quext{weakly in }\,\HUH,
\end{align}
for suitable limit functions $u$, $\chi$, $b$.
Then, the Aubin-Lions compactness Lemma and the usual
monotonicity argument \cite[Prop.~1.1, p.~42]{Ba}
permit to see that $b=b(\chi)$
a.e.~in~$(0,T)\times \Omega$. Moreover, the above 
relation suffice to pass to the limit in system
\eqref{phase1}-\eqref{phase2}.

Taking the limit in \eqref{calore} and in \eqref{identif} 
is a bit more involved. Actually, \eqref{st-energy} and
a comparison of terms in \eqref{calore} give
\begin{equation}\label{coa5}
  \teta_{n,t} \to \teta_t \quext{weakly in }\,\LDVp.
\end{equation}
Then, integrating in time and using the $V'$-convergence
in \eqref{tetaznn2}, we obtain more precisely
\begin{equation}\label{coa6}
  \teta_n \to \teta \quext{weakly in }\,\HUVp.
\end{equation}
This is sufficient to take the limit of equation
\eqref{calore}, but not of \eqref{identif}. 
Actually, to identify $\teta$ in terms of $u$, we 
have to work a little bit more. Namely, we 
have to integrate \eqref{calore} with respect 
to time both at the $n$-level and in the limit
and then test, respectively, by $u_n$ and by 
$u$. Notice that, even at the limit level,
the use of $u$ as a test function is guaranteed
by the fact that $u\in\LDV$ and all terms
in \eqref{calore} lie at least in $\LDVp$.

Then, at the $n$-level, we obtain
\begin{equation}\label{giro11}
  \iTT \duavb{\teta_n,u_n}
   + \frac12 \big\| \nabla (1*u_n) (t) \big\|^2
   = \iTT \duavb{\teta\zzn+\chi\zzn-\chi_n,u_n},
\end{equation}
where $*$ denotes convolution in time.
Taking the supremum limit in the above 
relation and comparing the result with the limit 
equation, we obtain
\begin{equation}\label{giro12}
  \limsup_{n\nearrow+\infty} \iTT \duavb{\teta_n,u_n}
   \le \iTT \duavb{\teta,u}.
\end{equation}
Since 
\begin{equation}\label{giro13}
  \teta_n = -\frac1{u_n}
   \in \devv J(u_n)\,\,\,\hbox{ a.e. in } (0,+\infty),
\end{equation}
(notice that we used here property \eqref{weaksubd2}),
relations \eqref{coa1} and \eqref{coa6} and 
the standard monotonicity argument \cite[Prop.~1.1, p.~42]{Ba},
applied here in the duality pairing between $V$ and $V'$,
permit to obtain \eqref{identifw}, which concludes the 
proof of existence in the conserved case.

\smallskip

\noindent%
{\bf Differences occurring in the non-conserved case.}~~%
At the level of estimates, the only relevant difference 
is in the energy relation, which is now obtained 
testing \eqref{calore} by $1+u$ and \eqref{phase} 
by $\chi_t$. This gives
\begin{equation}\label{energync}
  \ddt \calE + \| \nabla u \|^2
   + \| \chi_t \|^2 
   = 0.
\end{equation}
Thus, we have the $H$-norm of $\chi_t$ rather 
than the $V'$-norm on the \lhs\ (and consequently
we obtain \eqref{regochiw} in place 
of \eqref{regochi}). Estimates
\eqref{st-21}-\eqref{st-23} hold without variations
while \eqref{st-24} makes no longer sense.
Notice that \eqref{st-23} can now be obtained 
testing directly \eqref{phase} by $b(\chi)$.
The passage to the limit is analogous.

\smallskip

\noindent%
{\bf Proof of \eqref{energy-teo}.}~~%
We first observe that \eqref{identifw} is equivalent to
\begin{equation}\label{identifwstar}
  u \in \devvs J^*(\teta),
\end{equation}
almost everywhere in $(0,T)$. Then,
the standard integration by parts formula
\cite[p.~73]{Br}, applied in the duality between $V'$ and $V$, 
gives
\begin{equation}\label{abscon}
  J^*(\teta)\in AC([0,T]), \quad
  \duavb{\teta_t,u} = \ddt J^*(\teta).
\end{equation}
Thanks to this formula, in the non-conserved case
for any $V'$-solution we are allowed to test \eqref{calore}
by $u$ and \eqref{phase} by $\chi_t$. Integrating
over $(0,t)$ for arbitrary $t>0$, we obtain
exactly \eqref{energy-teo}. In the conserved case,
instead, we have to test \eqref{calore} by $u$, 
\eqref{phase1} by $w$ and \eqref{phase2} by $\chi_t$. 
Note that this is still possible for any
$V'$-solutions. Indeed, thanks to \eqref{regochiw} and 
the properties of $A$, we have that $w\in \LDV$. 
Thus, $w$ can be used as a test function in \eqref{phase1}
(which is a relation in $\LDVp$) and $\chi_t$ 
can be used as a test function in \eqref{phase2}
(which is a relation in $\LDV$ thanks to 
the above discussion). However, we have to notice 
that, while
\begin{equation}\label{separreg}
   A\chi + b(\chi) \in \LDV,
\end{equation}
it is not expected to be true that, separately,
$A\chi \in \LDV$ and $b(\chi)\in \LDV$. Nevertheless,
as shown, e.g., in \cite[Lemma~4.1]{RS}, property 
\eqref{separreg} is sufficient to prove that
\begin{equation}\label{separreg2}
  \duavb{ \chi_t , A\chi + b(\chi) }
   = \ddt \io \Big ( \frac12 |\nabla \chi|^2 + \bciapo(\chi) \Big)
\end{equation}
almost everywhere in $(0,T)$. Thus, we still
have \eqref{energy-teo}, of course with $\| \nabla w \|^2$
in place as $\| \chi_t \|^2$.

\smallskip

\noindent%
{\bf Proof of uniqueness.}~~%
It works exactly as in the standard case (so, we
just sketch it for the conserved model). Namely, we can
take a couple of solutions $(\teta_1,u_1,\chi_1)$,
$(\teta_2,u_2,\chi_2)$ starting from the same initial datum,
write the system for both solutions and take the difference.
Then, setting $(\teta,u,\chi):=(\teta_1,u_1,\chi_1)-(\teta_2,u_2,\chi_2)$,
we integrate (the difference of) \eqref{calore} in 
time and test it by $u$. Moreover, we test 
(the difference of) \eqref{phase1} by
$A^{-1}\chi$ (note that $\chi$ has zero-mean value, so
$A^{-1}$ is well-defined) and the difference of 
\eqref{phase2} by $\chi$. Collecting everything and
noting that two couples of terms cancel, we obtain
\begin{equation}\label{uniq11}
  \duav{\teta,u} + \frac12 \ddt \big( \| \nabla (1*u) \|^2
   + \| \chi \|_{V'}^2 \big)
   + \| \nabla \chi \|^2 
   \le \| \chi \|^2
   \le \frac12 \| \nabla \chi \|^2
   + c \| \chi \|_{V'}^2
\end{equation}
Noting that, by monotonicity, 
$\duav{\teta,u}=\duav{\teta_1-\teta_2,u_1-u_2}\ge 0$, 
the thesis follows then from Gronwall's lemma.

%%%%%%%%%%%%%%%%%%%%%%%%%%%%%%%%%%%%%%%%%%%%%%%%%%%%%%%%%%%%%%%%%%%

\subsection{Proof of Theorem~\ref{teo:asy:u}}

We start by deducing an additional a-priori estimate. 
As before, we present it just in the conserved case,
the variations in the non-conserved case being
given at the end.

\smallskip

\noindent%
{\bf Second estimate -- local version.}~~%
We test \eqref{calore} by $t u_t=t \teta_t/\teta^2$
and add the result to \eqref{phase1} multiplied by
$t w_t$. Then, we add also the time derivative of \eqref{phase2}
multiplied by $t \chi_t$. We obtain
\begin{align}\no
  & \ddt \Big( \frac{t}2 \| \nabla u \|^2
   + \frac{t}2 \| \nabla w \|^2 \Big)
   + t \io \frac{\teta_t^2}{\teta^2}
   + t \| \nabla \chi_t \|^2
   + t \io b'(\chi)\chi_t^2 \\
 \label{conto31} 
  & \mbox{}~~~~~
   \le \frac12\big( \| \nabla u \|^2
   + \| \nabla w \|^2 \big)
   + t \| \chi_t \|^2.
\end{align}
Then, noting that
\begin{equation}\label{conto32}
  t \| \chi_t \|^2
   \le \frac{t}2 \| \nabla\chi_t \|^2
   + c t \| \nabla w \|^2,
\end{equation}
integrating \eqref{conto31} between $0$ and $t\in(0,1]$
and taking advantage of \eqref{st-energy}, we obtain
\begin{equation}\label{st-31}
  t \| \nabla u(t) \|^2
   + t \| \nabla w(t) \|^2 
   + t \| \chi_t(t) \|_{V'}^2 
   + \int_0^t s \| \nabla \chi_t(s) \|^2\,\dis 
   \le c \EE_0.
\end{equation}

\smallskip

\noindent%
{\bf Second estimate -- global version.}~~%
We test \eqref{calore} by $u_t=\teta_t/\teta^2$
and add the result to \eqref{phase1} multiplied by
$w_t$. Then, we add also the time derivative of \eqref{phase2}
multiplied by $\chi_t$. Proceeding as above, we obtain
\begin{equation}\label{conto31b} 
  \ddt \Big( \frac{1}2 \| \nabla u \|^2
   + \frac{1}2 \| \nabla w \|^2 \Big)
   + \io \frac{\teta_t^2}{\teta^2}
   + \frac12 \| \nabla \chi_t \|^2
   + \io b'(\chi)\chi_t^2 
  \le c \| \nabla w \|^2.
\end{equation}
Integrating between $1$ and $t\ge 1$
and recalling \eqref{st-energy} and \eqref{st-31}, 
we infer  
\begin{equation}\label{st-31b}
  \| \nabla u \|_{L^\infty(1,\infty;H)}
   + \| \nabla w \|_{L^\infty(1,\infty;H)}
   + \| \chi_t \|_{L^\infty(1,\infty;V')}
   + \int_1^\infty \| \nabla \chi_t(t) \|^2\,\dit 
   \le Q(\EE_0).
\end{equation}
Using again the logarithmic Poincar\'e inequality
(Lemma~\ref{log-poincare}), we also have
\begin{equation}\label{st-31c}
  \| u \|_{L^\infty(1,\infty;V)}
   \le Q(\EE_0),
\end{equation}
i.e., the first \eqref{asym-u}.

In the non-conserved case, the procedure is similar. In place of
\eqref{st-31b} we rather obtain
\begin{equation}\label{st-31b-nc}
  \| \nabla u \|_{L^\infty(1,\infty;H)}
   + \| \chi_t \|_{L^\infty(1,\infty;H)}
   + \int_1^\infty \| \nabla \chi_t(t) \|^2\,\dit 
   \le Q(\EE_0).
\end{equation}
As a further consequence, we can look at equation \eqref{phase2}
in the conserved case (\eqref{phase} in the non-conserved case,
respectively). Thanks to estimates \eqref{st-31b}-\eqref{st-31c} 
for $u$ and $w$ (respectively, to estimate \eqref{st-31b-nc}
for $u$), applying standard regularity results for elliptic
equations with monotone nonlinearities, we then obtain 
\eqref{asym-chi}. 

\smallskip

\noindent%
{\bf Asymptotic uniform regularity of $u$.}~~%
Our aim is now to show the second \eqref{asym-u}. 
The key step is represented by the following lemma:
\bele\label{moser-u}
 Let $u$ be a solution of the problem
 \begin{equation}\label{eq-abstr}
   \teta_t + A u = f, \qquad \teta=-1/u
 \end{equation}
 over the time interval $(S,S+2)$, where we additionally
 assume that 
 \begin{equation}\label{inpiu}
   \| u \|_{L^3(S,S+2;L^{3/2}(\Omega))} \le M, \qquad
    \| f \|_{L^2(S,S+2;L^{3+\epsi}(\Omega))} \le F,
 \end{equation}
 for some (given) constants $M>0$, $F>0$ and some
 $\epsi>0$. Moreover, let us assume that 
 \begin{equation}\label{uSL1}
   u(S)\in L^1(\Omega).
 \end{equation}
 Then, we have
 \begin{equation}\label{st-abstr-u}
   \| u \|_{L^\infty((S+1,S+2)\times \Omega)} \le 
    Q\big( F, M, \| u(S) \|_1\big).
 \end{equation}
\enle
\begin{proof}
We test \eqref{eq-abstr} by $-|u|^{p+1}$, where $p\ge 1$ 
will be specified later (although $-|u|^{p+1}$ needs not
necessarily be an admissible test function, the procedure
could be easily justified by truncation arguments, we
omit the details). This gives
\begin{equation}\label{lemma1-1}
  \frac1p \ddt \| u \|_p^p
   + \frac{4(p+1)}{(p+2)^2} \big\| \nabla |u|^{\frac{p+2}2} \big\|^2
   \le \io |f| |u|^{p+1}.
\end{equation}
We then set $r:=\frac{3+\epsi}{2+\epsi}$ to be the conjugate 
exponent of $3+\epsi$. Then, multiplying by $p$, 
we can estimate the right hand side as
\begin{equation}\label{lemma1-2}
  p\io |f| | u |^{p+1}
   \le p \|f\|_{3+\epsi} \big\| |u|^{p+1} \big\|_r
   = p \| f \|_{3+\epsi} \| u \|_{r(p+1)}^{p+1}.
\end{equation}
Then, in order to recover the full $V$-norm from
the gradient term, we add 
\begin{equation}\label{lemma1-3}
  \big\| |u|^{\frac{p+2}2} \big\|_1^2
   = \| u \|_{\frac{p+2}2}^{p+2}
\end{equation}
to both hands sides of \eqref{lemma1-1}.
Integrating \eqref{lemma1-1} over $(\tau,t)$,
for $t$ a generic point in $(\tau,S+2)$
and choosing, for the first iteration, $p=1$ 
and $\tau=S$, we obtain
\begin{align} \no
  & \| u \|_{L^{\infty}(\tau,S+2;L^{p}(\Omega))}^p
   + \| u \|_{L^{p+2}(\tau,S+2;L^{3p+6}(\Omega))}^{p+2} \\
 \no 
  & \mbox{}~~~~~~~~~~
   \le c \| u(\tau) \|_{L^p(\Omega)}^p
   + c p \| f \|_{L^2(S,S+2;L^{3+\epsi}(\Omega))}
    \| u \|_{L^{2(p+1)}(\tau,S+2;L^{r(p+1)}(\Omega))}^{p+1}
   + c\int_\tau^{S+2} \| u(s) \|_{\frac{p+2}2}^{p+2} \,\dis \\
 \label{lemma1-4}
 & \mbox{}~~~~~~~~~~
   \le c p F
    \| u \|_{L^{2(p+1)}(\tau,S+2;L^{r(p+1)}(\Omega))}^{p+1}
    + Q\big( M, \| u(S) \|_p\big),
\end{align}
where in the last inequality we took advantage of 
\eqref{inpiu} using that $p=1$ and $\tau=S$.

Being non-restrictive to assume that $u\ge 1$ almost
everywhere (otherwise, we can replace $u$ with 
$\max\{u,1\}$), we can then define 
\begin{equation}\label{Jp}
  J_p^p := \| u \|^{p}_{L^\infty(\tau_p,S+2;L^p(\Omega))} 
   + \| u \|^{p}_{L^{p+2}(\tau_p,S+2,L^{3p+6}(\Omega))}, 
\end{equation} 
where, for now, we take $\tau_p=\tau=0$. Then,
by interpolation we obtain
\begin{equation}\label{nota-11}
  \| u \|_{L^{2(q+1)}(\tau_p,S+2;L^{r(q+1)}(\Omega))}
   \le \| u \|_{L^\infty(\tau_p,S+2;L^p(\Omega))}^{\alpha}
    \| u \|_{L^{p+2}(\tau_p,S+2,L^{3p+6}(\Omega))}^{1-\alpha},
\end{equation}
for some $\alpha\in(0,1)$.
Then, raising to the power $p$ and using the Young inequality 
with exponents $P=1/\alpha$ e $Q=1/(1-\alpha)$, we get
\begin{equation}\label{nota-12}
  \| u \|_{L^{2(q+1)}(\tau_p,S+2;L^{r(q+1)}(\Omega))}^p
   \le \alpha \| u \|_{L^\infty(\tau_p,S+2;L^p(\Omega))}^p
    + (1-\alpha) \| u \|_{L^{p+2}(\tau_p,S+2,L^{3p+6}(\Omega))}^p,
\end{equation}
which implies, upon dividing by $\max\{\alpha,1-\alpha\}$ 
(that is different from $0$ and $1$)
\begin{equation}\label{Jpubis}
   J_p^p \ge \| u \|_{L^{2(q+1)}(\tau_p,S+2;L^{r(q+1)}(\Omega))}^p,
\end{equation}
where the index $q$ and the interpolation exponent $\rho$ 
are given by the system
\begin{equation} \label{interpol12}
  \displaystyle
   \begin{cases}
     \frac{1-\rho}{p+2} = \frac{1}{2(q+1)},\\
      \frac{\rho}{p} + \frac{1-\rho}{3(p+2)}\,=\,\frac{1}{r(q+1)}.
   \end{cases}
\end{equation}
Dividing the second equation in \eqref{interpol12} 
by the first one, we actually have
\begin{equation} \label{interpol13}
  \Big(\frac{\rho}{p} + \frac{1-\rho}{3(p+2)} \Big)\frac{p+2}{1-\rho}=\frac{2}{r},
\end{equation}
whence 
\begin{equation} \label{interpol14}
  \frac{\rho}{p}\frac{p+2}{1-\rho}=\frac{2}{r}-\frac{1}{3}=:K_\epsi,
\end{equation}
and it is easy to compute
\begin{equation} \label{interpol15}
  K_\epsi=\frac{9+7\epsi}{9+3\epsi}.    
\end{equation} 
From \eqref{interpol14} and the first equation in \eqref{interpol12}, we then 
have
\begin{equation} \label{interpol16}
  \rho = \frac{ \frac{p}{p+2} K\ee }{ 1 + \frac{p}{p+2} K\ee } 
   \in (0,1) \quad \perogni p\ge 1.
\end{equation}
Being 
\begin{equation} \label{interpol17}
  1 - \rho = \frac{ 1 }{ 1 + \frac{p}{p+2} K\ee },
\end{equation}
we then obtain fron the first \eqref{interpol12}
\begin{align} \no
  q & = \frac12 \frac{p+2}{1-\rho} - 1
    = \frac12 \Big(1 + \frac{p}{p+2} K\ee \Big) (p+2) - 1\\
  \label{interpol18}
  & = \frac{ K\ee + 1 }2 p = \frac{9+5\epsi}{9+3\epsi} p =: H p,
\end{align}
where, obviously, $H=H(\epsi)>1$ whenever $\epsi>0$. 

\smallskip

Given that $p_0=1$, let us set, by induction, 
$p_{i+1}=H p_i = H^{i+1}$.
Then, let $i\ge 0$ and let us rewrite \eqref{lemma1-4}
by taking $p=p_{i+1}$ and $\tau=\tau_{i+1}$
(the latter will be chosen below).
Setting also, for brevity, $J_i:=J_{p_i}$, we
then obtain, thanks also to \eqref{Jpubis},
\begin{equation} \label{moser10}
  J_{i+1}^{p_{i+1}} 
    \le c \| u(\tau_{i+1}) \|_{p_{i+1}}^{p_{i+1}}
    + c p_{i+1} F J_i^{p_{i+1}+1} 
   + c\int_{\tau_{i+1}}^2 \| u(s) \|_{\frac{p_{i+1}+2}2}^{p_{i+1}+2} \,\dis.
\end{equation}
Now, let (for instance), for $i\ge 1$,
\begin{equation} \label{sigmai}
  \sigma_i=\frac6{\pi^2i^2}, \quext{so that }\,
   \sum_{i=1}^\infty \sigma_i = 1.
\end{equation}
Then, we observe that, given $\tau_i$, we can choose 
$\tau_{i+1} \in (\tau_i,\tau_i+\sigma_{i+1})$
such that
\begin{align} \no
  \| u(\tau_{i+1}) \|_{p_{i+1}}^{p_{i}}
  & \le \frac1{\sigma_{i+1}}\int_{\tau_i}^{\tau_i+\sigma_{i+1}}
   \| u(t) \|_{p_{i+1}}^{p_{i}}\,\dit\\
  \label{moser10b}
  & \le c i^2 \int_{\tau_i}^{\tau_i+\sigma_{i+1}}
   \| u(t) \|_{3(p_i+2)}^{p_i}\,\dit
  \le c i^2 J_i^{p_i},
\end{align}
where we used that $p_{i+1}\le 3(p_i+2)$.

Analogously, we have that
\begin{equation} \label{moser10c}
  c\int_{\tau_{i+1}}^{S+2} \| u(s) \|_{\frac{p_{i+1}+2}2}^{p_{i+1}+2} \,\dis
   \le 2c \| u \|_{L^\infty(\tau_i,S+2;L^{p_i}(\Omega))}^{p_{i+1}+2}
   \le c J_i^{p_{i+1}+2}.
\end{equation}
Collecting \eqref{moser10b}-\eqref{moser10c},
\eqref{moser10} gives
\begin{equation} \label{moser11}
  J_{i+1}^{p_{i+1}} 
    \le c i^{2H} J_i^{p_{i+1}}
    + c p_{i+1} F J_i^{p_{i+1}+1} 
    + c J_i^{p_{i+1}+2}.
\end{equation}
Thus, 
%noting that, as $i$ becomes large, $C_{\rho_i}$ does not 
%degenerate (cf.~\eqref{interpol16}), 
we finally obtain
\begin{equation} \label{moser12}
  J_{i+1}^{H^{i+1}} \le c \big( i^{2H} + H^{i+1} F + 1 \big)
   J_i^{H^{i+1}+2}.
\end{equation}
Thus, setting $\eta_i:=(H^i+2)/H^i$, we have
\begin{equation} \label{moser13}
  J_{i+1} \le B_i^{H^{-(i+1)}}
   J_i^{\eta_{i+1}}, \quext{where }\,
   B_i:= c \big( i^{2H} + H^{i+1} F + 1 \big)
\end{equation}
whence a simple induction argument (see, e.g., \cite{Sc})
permits to obtain \eqref{st-abstr-u}.
\end{proof}
\noindent
{\bf Conclusion of the proof.}~~%
To obtain the uniform boundedness of $u$ it is now 
sufficient to notice that, by the first \eqref{asym-u},
$u(t)\in L^6(\Omega)$ for all $t\ge 1$. Then, \eqref{uSL1} holds.
Moreover, \eqref{inpiu} are a consequence of \eqref{st-31b}
(or, in the non-conserved case, \eqref{st-31b-nc}),
which gives the required regularity for $f=-\chi_t$,
and of \eqref{st-31c}, which gives the required regularity for $u$.
The second \eqref{asym-u} is then a consequence of the lemma
(applied with the choice of $S=2$). It is also worth noting that, 
at the level of $V'$-solution, the identification property in
\eqref{eq-abstr} needs not hold in the strong (pointwise)
form (but just in the ``weak'' sense \eqref{identifw}). 
However, one can apply Lemma~\ref{moser-u} at the 
$n$-regularized level, and then pass to the limit
noting that the procedure yields $n$-uniform estimates.

To complete the proof we have to show \eqref{separ} in the 
non-conserved case. Of course, such a property is significant
only in the case when $I$, the domain of $b$, does not coincide
with the real line, i.e., we are in presence of a singular 
potential (like the logarithmic one \eqref{potlog}). 
Otherwise, \eqref{separ} is (also in the conserved case, of course)
an immediate consequence of \eqref{asym-chi}. 
 
That said, let us prove the upper bound, the lower one working in
a similar way. Being, by the second \eqref{asym-u}, $|u(t,x)|\le U$ 
for some $U=Q(\EE_0)$ and a.e.~$(t,x)\in(2,\infty)\times \Omega$, 
we can then apply the comparison principle to \eqref{phase}.
This gives that $\chi$ is bounded from above by the 
solution $\chi^+$ to the forward Cauchy problem
\begin{equation} \label{phasecomp}
  \chi^+_t + b(\chi^+) - \chi^+ = U,
   \quad \chi^+(2)= I^+,
   \quext{where }\,I^+:=\sup I.
\end{equation}
Actually, by the last \eqref{hpb}, 
$\lim_{r\nearrow \sup I} b(r) - r - U= +\infty$.
Thus, $\chi^+(t) \le I^+ - \delta$ for all $t\ge 3$ 
and some $\delta > 0$. Moreover, this bound is uniform
in time. Then, \eqref{separ} is a consequence of the 
comparison principle. This concludes the proof 
of the theorem.

%%%%%%%%%%%%%%%%%%%%%%%%%%%%%%%%%%%%%%%%%%%%%%%%%%%%%%%%%%%%%%%%%%%%%%%%%%%%%%%%%%%%
%%%%%%%%%%%%%%%%%%%%%%%%%%%%%%%%%%%%%%%%%%%%%%%%%%%%%%%%%%%%%%%%%%%%%%%%%%%%%%%%%%%%

\subsection{Proof of Theorem~\ref{teo:esi:w2}} 
\label{sec:vplu}

We start by giving the proof in the 
non-conserved case, the variations occurring in
the conserved case being outlined at the end.
Then, we know that any $V'$-solution 
$(\teta,u,\chi)$ satisfies
the energy equality \eqref{energy-teo}.
Analogously, if $(\teta_n,u_n,\chi_n)$ is the 
approximating solution constructed in the existence
proof, the analogue of \eqref{energy-teo}
reads
\begin{align}\no
  & \itt \big( \| \nabla u_n \|^2
   + \| \chi_{n,t} \|^2 \big)
   = J^*(\teta\zzn) - J^*(\teta_n(t)) 
   + \duav{\teta\zzn,1} - \duav{\teta_n(t),1}\\
 \label{energync-n}  
   & \mbox{}~~~~~
   + \io\Big(\chi_0 + \frac12|\nabla\chi_0|^2 + \bhat(\chi_0)- \frac12\chi_0^2 \Big)
   - \io\Big(\chi_n(t) + \frac12|\nabla\chi_n(t)|^2 + \bhat(\chi_n(t))- \frac12\chi_n(t)^2 \Big).
\end{align}
Our task is now to compute the supremum limit of 
\eqref{energync-n} and compare it with \eqref{energy-teo}.
Then, we firstly observe that, by \eqref{coa2}, 
the Aubin-Lions lemma, and lower
semicontinuity of $\bhat$,
\begin{equation}\label{limsup11}
  \io\Big(\chi(t) + \frac12|\nabla\chi(t)|^2 + \bhat(\chi(t))- \frac12\chi(t)^2 \Big)
  \le \liminf_{n\nearrow\infty}\io\Big(\chi_n(t) + \frac12|\nabla\chi_n(t)|^2 
    + \bhat(\chi_n(t))- \frac12\chi_n(t)^2 \Big).
\end{equation}
Analogously, using convexity and lower semicontinuity of the functional $J^*$
w.r.t.~the $V'$-norm, \eqref{tetaznn2}, and the fact that, by \eqref{coa6},
$\teta_n(t)$ tends to $\teta(t)$ weakly in $V'$ for {\sl all}\/ $t\in[0,T]$,
we obtain
\begin{equation}\label{limsup12}
  J^*(\teta(t)) - \duav{\teta_0,1} + \duav{\teta(t),1}
   \le \liminf_{n\nearrow\infty} \Big( J^*(\teta_n(t)) 
   - \duav{\teta\zzn,1} + \duav{\teta_n(t),1} \Big).
\end{equation}
So, it remains to prove that 
\begin{equation}\label{limsup13}
  \limsup_{n\nearrow\infty} J^*(\teta\zzn) 
   \le J^*(\teta_0).
\end{equation}
The proof of this fact is actually
a bit more involved. We prepare a Lemma
\bele\label{lemma:sergey}
 Let $v\in V'\cap L^1(\Omega)$ such that $j^*(v) \in L^1(\Omega)$,
 $j$ and $j^*$ being given by\/ \eqref{defij}. Then,
 \begin{equation}\label{ls-11}
   J^*(v) = \io j^*(v).
 \end{equation}
\enle
\begin{proof}
Let $z$ belong to the domain of $J$, namely let $z\in V$ with 
$J(z)<+\infty$. Then, for a.e.~$x\in \Omega$, by
definition of subdifferential in $\RR$, we have
\begin{equation}\label{ls-12}
  z(x) v(x) \le j(z(x)) + j^*(v(x)).
\end{equation}
Integrating over $\Omega$, we would {\sl formally} get
\begin{equation}\label{ls-13}
  \io zv \le J(z) + \io j^*(v).
\end{equation}
However, the integral on the \lhs\ could make no sense
since the function $zv$ could not belong to $L^1(\Omega)$. 
Nevertheless, it is simple (see, e.g., \cite[Lemma~2.2]{BCGG} and
\cite[Lemma~2.1]{Ba}) to see that, in place of \eqref{ls-13},
there holds
\begin{equation}\label{ls-14}
  \duav{ z, v } \le J(z) + \io j^*(v).
\end{equation}
Passing to the supremum w.r.t.~$z$ varying in the 
domain of $J$, we then get the $\le$ sign in \eqref{ls-11}.

To prove the converse, we first let, for 
$v$ as in the statement and $N\in(1,\infty)$,
\begin{equation}\label{ls-21}
  z_N:= - \max\big\{N^{-1},\min\{v^{-1},N\}\big\}.
\end{equation}
Then, for $\epsilon\in (0,1)$, we regularize by
singular perturbation as in \eqref{tetaznn}:
\begin{equation}\label{ls-22}
  z\eeN + \epsilon A z\eeN = z_N.
\end{equation}
Then,
\begin{align}\no
  J^*(v) & \ge \duav{v,z\eeN} - J(z\eeN)
   = \io v z\eeN - J(z\eeN)
   = \io v z\eeN + \io \log( - z\eeN)\\
 \no
  & \ge \io v z\eeN + \io \log ( - z_N)
   \to^{\epsilon\searrow0} \io v z_N + \io \log ( - z_N)\\
 \label{ls-23}
  & \to^{N\nearrow\infty} \io ( -1 + \log v^{-1} )
    = \io j^*(v),
\end{align}
where the first inequality follows from definition
of conjugate function, the second equality from
the fact that $z\eeN$ is smooth and bounded, the third 
equality is trivial, the fourth inequality comes
from \eqref{ls-22}, the convergence $\epsi\searrow0$
from standard properties of elliptic systems,
and the convergence $N\nearrow\infty$ from 
Lebesgue's theorem. This proves the $\ge$ in
\eqref{ls-11} and the lemma.
\end{proof}
\noindent%
Then, by definition of  conjugate function 
(recall \eqref{defij}),
\begin{align}\no
  J^*(\teta\zzn) 
   & = \sup_{v\in V} \duavb{\teta\znn,v} - J(v)
   \le \sup_{v\in H} \big(\teta\znn,v\big) - J(v)\\
  \label{limsup14}
  & = \io -1 - \log\teta\znn
% \io - \log \teta\zzn 
  \le \io - 1 - \log \teta_0,
\end{align}
where the last inequality follows from
\eqref{tetaznn3}. In particular, we have that
$-1 -\log\teta_0\in L^1(\Omega)$.
Thus, applying the above Lemma, we obtain
\begin{equation}\label{Jsvp3}
  J^*(\teta_0) = \io (-1 - \log\teta_0),
\end{equation}
whence, computing the supremum limit of \eqref{limsup14}, 
\eqref{limsup13} follows. Thus, we finally end up with
\begin{equation}\label{giro21}
  \limsup_{n\nearrow\infty} 
   \itt \big( \| \nabla u_n \|^2
   + \| \chi_{n,t} \|^2 \big)
  \le \itt \big( \| \nabla u \|^2
   + \| \chi_{t} \|^2 \big),
\end{equation}
whence, recalling \eqref{coa1}-\eqref{coa4}, we get in particular
\begin{equation}\label{giro22}
  \nabla u_n,~\chi_{n,t}
   \to \nabla u,~\chi_{t}
   \quext{strongly in }\,\LDH.
\end{equation}
As a final step of our procedure, we shall prove that
$\{\teta_n\}$ is a Cauchy sequence in $C^0([0,T];L^1(\Omega))$.
Actually, writing equation \eqref{calore} for a couple
of indexes $n$ and $m$ and taking the difference, we 
obtain
\begin{equation}\label{calorediff}
  \teta_{n,t} - \teta_{m,t} 
   + A (u_n - u_m) 
   = \chi_{n,t} - \chi_{m,t}.
\end{equation}
Thus, testing by $\sign(\teta_n-\teta_m)$,
noticing that, by monotonicity, 
$\sign(\teta_n-\teta_m)=\sign(u_n-u_m)$,
and applying the Brezis-Strauss theorem
\cite[Lemma~2]{BS}, we arrive at 
\begin{equation}\label{BS11}
  \ddt \| \teta_n - \teta_m \|_{L^1(\Omega)}
   \le \| \chi_{n,t} - \chi_{m,t} \|_{L^1(\Omega)},
\end{equation}
whence, integrating in time and using the 
strong convergences \eqref{tetaznn2} and 
\eqref{giro22}, we end up with
\begin{equation}\label{BS12}
  \teta_n \to \teta
   \quext{strongly in }\,C^0([0,T];L^1(\Omega)).
\end{equation}
Thus, in particular, we have that, for all $t\in[0,T]$,
$\teta(t)\in L^1(\Omega) \cap \devv J(u(t))$, whence 
the pointwise identification \eqref{identif} follows 
from \eqref{weaksubd3}. This concludes the proof in the 
non-conserved case.

\smallskip

\noindent%
{\bf Conserved case.}~~%
To conclude the proof, we outline the differences occurring
in the conserved case, which only regard the above $L^1$-argument.
Actually, the convergence of $\chi_{n,t}$ in 
\eqref{giro22} is now replaced by
\begin{equation}\label{giro22w}
  \nabla w_n \to \nabla w,
   \quext{strongly in }\,\LDH.
\end{equation}
Of course, thanks to the properties of $A$ this also
gives
\begin{equation}\label{giro23w}
  \chi_{n,t} \to \chi_t
   \quext{strongly in }\,\LDVp,
\end{equation}
which, however, is not sufficient to proceed as before.
On the other hand, we can rely on estimates 
\eqref{st-31} and \eqref{st-31b-nc} which tell us
that
\begin{equation}\label{serg11}
   \big\| t^{1/2} \chi_{n,t} \big \|_{\LDV}
   + \big\| t^{1/2} \chi_{n,t} \big\|_{\LIVp} \le c.
\end{equation}
Thus, by interpolation, 
\begin{equation}\label{serg12}
  \big\| t^{1/2} \chi_{n,t} \big\|_{L^{\frac4{2-\epsilon}}(0,T;H^{1-\epsilon}(\Omega))} \le c
\end{equation}
for all $\epsilon\in(0,1)$. Coming back to \eqref{BS11}, we now have that,
for all $t\in[0,T]$,
\begin{align}\no
  \| \chi_{n,t} - \chi_{m,t} \|_{L^1(\Omega)}
   & \le c \| \chi_{n,t} - \chi_{m,t} \|
    \le c \| \chi_{n,t} - \chi_{m,t} \|_{V'}^{\frac{1-\epsilon}{2-\epsilon}}
      \| \chi_{n,t} - \chi_{m,t} \|_{H^{1-\epsilon}(\Omega)}^{\frac1{2-\epsilon}}\\
\label{serg13}      
   & \le c t^{-{\frac1{2(2-\epsilon)}}} \| \chi_{n,t} - \chi_{m,t} \|_{V'}^{\frac{1-\epsilon}{2-\epsilon}} \cdot
     t^{\frac1{2(2-\epsilon)}} \| \chi_{n,t} - \chi_{m,t} \|_{H^{1-\epsilon}(\Omega)}^{\frac1{2-\epsilon}}.
\end{align}
Then, integrating over $(0,T)$, we arrive at
\begin{align}\no
  \| \chi_{n,t} - \chi_{m,t} \|_{L^1((0,T)\times\Omega)}
   & \le c \Big\| t^{-{\frac1{2(2-\epsilon)}}} \| \chi_{n,t} - \chi_{m,t} \|_{V'}^{\frac{1-\epsilon}{2-\epsilon}} 
       \Big\|_{L^{\frac43}(0,T)} 
     \Big\| t^{\frac1{2(2-\epsilon)}} \| \chi_{n,t} - \chi_{m,t} \|_{H^{1-\epsilon}(\Omega)}^{\frac1{2-\epsilon}} 
      \Big\|_{L^4(0,T)}\\
 \no
   & \le c \bigg( \iTT t^{-\frac2{3(2-\epsilon)}} 
      \| \chi_{n,t} - \chi_{m,t} \|_{V'}^{\frac{4(1-\epsilon)}{3(2-\epsilon)}} \bigg)^{\frac34}\\
 \no
   & \le c \big\| t^{-\frac2{3(2-\epsilon)}} \big\|_{L^2(0,T)}^{\frac34} 
     \Big\| \| \chi_{n,t} - \chi_{m,t} \|_{V'}^{\frac{4(1-\epsilon)}{3(2-\epsilon)}} \Big\|_{L^2(0,T)}^{\frac34}\\
 \label{serg14}  
   & \le c \big\| \chi_{n,t} - \chi_{m,t} \|_{L^{\frac{8(1-\epsilon)}{3(2-\epsilon)}}(0,T;V')}%
     ^{\frac{1-\epsilon}{2-\epsilon}},
\end{align}
where the second inequality is a consequence of \eqref{serg12},
the thirds follows from H\"older's inequality, and the fourth
holds provided that we take $\epsilon$ so small that $4<3(2-\epsilon)$.
In particular, using \eqref{giro23w} (note that $8(1-\epsilon)/3(2-\epsilon)$
is smaller than $2$ for $\epsilon$ as above), we obtain that
the \rhs\ tends to $0$ for large
$m$ and $n$. At this point the proof
goes on like in the non-conserved case.

%%%%%%%%%%%%%%%%%%%%%%%%%%%%%%%%%%%%%%%%%%%%%%%%%%%%%%%%%%%%%%%%%%%

\subsection{Proof of Theorem~\ref{teo:asy:teta}}

To start, we need to prove some further a-priori estimates 
holding under the additional assumption~\eqref{tetazl3}.
In particular, the key step will be that of showing
that the $L^{3+\epsi}$ regularity of $\teta$ is conserved
uniformly in time. To show this, we will use 
$L^p$-techniques in equation~\eqref{calore} (actually,
with $p=3+\epsi$). However, doing this will require
some care since, due to the low regularity of initial data,
the ``forcing term'' $\chi_t$ needs
not belong to $L^{3+\epsi}(\Omega)$ 
for small values of the time variable. However, we will
see that the $L^{3+\epsi}$-norm of $\chi_t(t)$ explodes,
as $t\searrow 0$, in a way which is sufficiently slow for
our purpose. As before, the proof is detailed just in the
conserved case. That said, we start with the

\smallskip

\noindent%
{\bf Third estimate -- local version.}~~%
To start, we test \eqref{calore} by $\teta^{2+\epsi}$,
to obtain
\begin{equation}\label{conto41}
  \ddt \| \teta \|_{3+\epsi}^{3+\epsi}
   \le c \io | \chi_t | \teta^{2+\epsi}
   \le c \| \chi_t \|_{3+\epsi}
   \| \teta \|_{3+\epsi}^{2+\epsi},
\end{equation}   
whence, clearly,
\begin{equation}\label{conto42}
  \ddt \| \teta \|_{3+\epsi}
   \le c \| \chi_t \|_{3+\epsi}.
\end{equation}   
Thus, using that
\begin{equation}\label{conto43}
  H^{\frac{3+3\epsi}{6+2\epsi}}(\Omega)
   \subset L^{3+\epsi}(\Omega)
\end{equation}   
and that
\begin{equation}\label{conto44}
  \frac{3+3\epsi}{6+2\epsi}
   = \alpha \times 1 + (1-\alpha)\times (-1),
   \quext{with }\,\alpha = \frac{9+5\epsi}{12+4\epsi},  
\end{equation}   
we obtain
\begin{align}\no
  \ddt \| \teta \|_{3+\epsi}
   & \le c \| \chi_t \|_{V}^{\frac{9+5\epsi}{12+4\epsi}}
    \| \chi_t \|_{V'}^{\frac{3-\epsi}{12+4\epsi}} \\
 \no  
   & \le c \big( t^{\frac{9+5\epsi}{24+8\epsi}}
             \| \chi_t \|_{V}^\frac{9+5\epsi}{12+4\epsi} \big)
    \big( t^{\frac{3-\epsi}{24+8\epsi}}
              \| \chi_t \|_{V'}^\frac{3-\epsi}{12+4\epsi} \big)
    t^{-\frac12} \\
 \no  
   & \le Q(\EE_0) \big( t^{\frac{9+5\epsi}{24+8\epsi}}
             \| \chi_t \|_{V}^\frac{9+5\epsi}{12+4\epsi} \big)
      t^{-\frac12} \\
 \label{conto45}
  & \le Q(\EE_0) \big( t \| \chi_t \|_{V}^2 
    + t^{-\frac{12+4\epsi}{15+3\epsi}} \big),
\end{align}
where \eqref{st-31} has also been exploited.
Note that, for $\epsi\in(0,1)$ the latter exponent lies in 
$(-1,0)$. Notice also that in the non-conserved case
the exponents are even better since it is sufficient
to interpolate between $V$ and $H$ (rather
than between $V$ and $V'$).
Thus, integrating \eqref{conto45} between $0$ 
and $1$, and using once more \eqref{st-31}, we infer
\begin{equation}\label{st-41}
  \| \teta \|_{L^\infty(0,1;L^{3+\epsi}(\Omega))} 
   \le Q\big(\EE_0,\|\teta_0\|_{3+\epsi}\big).
\end{equation}

\smallskip

\noindent%
{\bf Third estimate -- global version.}~~%
As before, we test \eqref{calore} by $\teta^{2+\epsi}$.
Taking now care also of the gradient term, we get
\begin{equation}\label{conto51}
  \ddt \| \teta \|_{3+\epsi}^{3+\epsi}
   + \kappa \big\| \nabla \teta^{\frac{1+\epsi}2} \big\|^2
   \le c \io | \chi_t | \teta^{2+\epsi}.
\end{equation}   
Adding also the inequality (which is true thanks to
\eqref{st-energy})
\begin{equation}\label{conto52}
  \kappa \big\| \teta^{\frac{1+\epsi}2} \big\|_1^2
   \le c \big(1 + \| \teta \|_1^2 \big) 
   \le Q(\EE_0),
\end{equation}   
we then get 
\begin{align}\no
  \ddt \| \teta \|_{3+\epsi}^{3+\epsi}
   + \kappa \big\| \teta^{\frac{1+\epsi}2} \big\|_6^2
   & \le c \| \chi_t \|_6 \big\| \teta^{\frac{3+\epsi}2} \big\|_{3/2}
      \big\| \teta^{\frac{1+\epsi}2} \big\|_6 + Q(\EE_0)\\
 \no     
   & \le c \| \chi_t \|_V^2\big\| \teta^{\frac{3+\epsi}2} \big\|_{3/2}^2
       + \frac\kappa2 \big\| \teta^{\frac{1+\epsi}2} \big\|_6^2 + Q(\EE_0)\\
 \no
   & \le c \| \chi_t \|_V^2 \big\| \teta \big\|_{\frac{9+3\epsi}4}^{3+\epsi}
       + \frac\kappa2 \big\| \teta^{\frac{1+\epsi}2} \big\|_6^2 + Q(\EE_0)\\
 \label{conto53}
   & \le c \| \chi_t \|_V^2 \big\| \teta \big\|_{3+\epsi}^{3+\epsi}
       + \frac\kappa2 \big\| \teta^{\frac{1+\epsi}2} \big\|_6^2 + Q(\EE_0).
\end{align}   
Thus, setting $Y(t):=\| \teta(t) \|_{3+\epsi}^2$,
%
% and noting
%that it is not restrictive to assume $Y(t)\ge 1$, 
%
we obtain from \eqref{conto53}
\begin{equation}\label{conto54}
  Y' + \kappa 
   \le c \| \chi_t \|_V^2 Y 
   + Q Y^{-\frac{1+\epsi}2},
\end{equation}   
where we wrote $Q$ in place of $Q(\EE_0)$.
Now, let us set
\begin{equation}\label{conto55}
  Z(t):=\max\Big\{ Y(t), \Big(\frac{Q}{\kappa}\Big)^{\frac{2}{1+\epsi}} \Big\},
\end{equation}   
so that it is clear that $Z$ satisfies, 
\begin{equation}\label{conto56}
  Z'(t) \le c \| \chi_t(t) \|_V^2 Z(t), 
  \quext{for a.e.~}\,t\ge 1, \qquext{and }\,
   Z(1) \le Q\big(\EE_0,\|\teta_0\|_{3+\epsi}\big).
\end{equation}   
Hence, integrating between $1$ and a generic $t\ge 1$ and 
recalling \eqref{st-31b}, we arrive at 
\begin{equation}\label{st-51}
  \| \teta \|_{L^\infty(1,\infty;L^{3+\epsi}(\Omega))}
   \le Q\big(\EE_0,\|\teta_0\|_{3+\epsi}\big),
\end{equation}   
i.e., the global analogue of \eqref{st-41}.

\smallskip

\noindent%
{\bf Asymptotic uniform regularity of $\teta$.}~~%
The key step is represented by the following counterpart
of Lemma~\ref{moser-u}:
\bele\label{moser-teta}
 Let $u$ be a solution of the problem
 \begin{equation}\label{eq-abstr2}
   \teta_t + A u = f, \qquad \teta=-1/u,
 \end{equation}
 over the time interval $(S,S+2)$, where we additionally
 assume that 
 \begin{equation}\label{inpiu-2}
   \| \teta \|_{L^\infty(S,S+2;L^{1}(\Omega))} \le M, \qquad
    \| f \|_{L^2(S,S+2;L^{3+\epsi}(\Omega))} \le F,
 \end{equation}
 for some (given) constants $M>0$, $F>0$ and some
 $\epsi>0$. Moreover, let us assume that 
 $\teta(S)\in L^{3+\epsi}(\Omega)$ for some $\epsi>0$.
 Then,
 \begin{equation}\label{st-abstr-teta}
   \| \teta \|_{L^\infty((S+1,S+2)\times \Omega)} \le 
    Q\big( F, M, \| \teta(S) \|_{3+\epsi} \big).
 \end{equation}
\enle
\begin{proof}
We test \eqref{eq-abstr} by $\teta^{p-1}$, where $p>3$ 
will be specified later. This gives
\begin{equation}\label{lemma2-1}
  \frac1p \ddt \| \teta \|_p^p
   + \frac{4(p-1)}{(p-2)^2} \big\| \nabla \teta^{\frac{p-2}2} \big\|^2
   \le \io |f| \teta^{p-1}.
\end{equation}
Setting as before $r:=\frac{3+\epsi}{2+\epsi}$ and multiplying by $p$, 
we have
\begin{equation}\label{lemma2-2}
  p\io |f| \teta^{p-1}
   \le p \|f\|_{3+\epsi}\| \teta^{p-1} \|_r
   = p \| f \|_{3+\epsi} \| \teta \|_{r(p-1)}^{p-1}
\end{equation}
and, in order to recover the full $V$-norm from
the gradient term, we add 
\begin{equation}\label{lemma2-3}
  \big\| \teta^{\frac{p-2}2} \big\|_1^2
   = \| \teta \|_{\max\left\{\frac{p-2}2,1\right\}}^{p-2}.
\end{equation}
Integrating \eqref{lemma2-1} over $(\tau,t)$,
for $t$ a generic point in $(\tau,S+2)$
and choosing, for the first iteration, $p=3+\epsi$ 
and $\tau=S$, we then have the analogue
of~\eqref{lemma1-4}:
\begin{align} \no
  & \| \teta \|_{L^{\infty}(\tau,S+2;L^{p}(\Omega))}^p
   + \| \teta \|_{L^{p-2}(\tau,S+2;L^{3p-6}(\Omega))}^{p+2} \\
 \no 
  & \mbox{}~~~~~~~~~~
   \le c \| u(\tau) \|_{L^p(\Omega)}^p
   + c p F
    \| \teta \|_{L^{2(p-1)}(\tau,S+2;L^{r(p-1)}(\Omega))}^{p+1}
   + c\int_\tau^{S+2} \| \teta(s) \|_{\max\left\{\frac{p-2}2,1\right\}}^{p-2} \,\dis \\
 \label{lemma2-4}
 & \mbox{}~~~~~~~~~~
   \le c p F
    \| \teta \|_{L^{2(p-1)}(\tau,t;L^{r(p-1)}(\Omega))}^{p-1}
    + Q\big( M, \| u_0 \|_1\big),
\end{align}
Now, the iteration scheme goes through similarly as before. 
Actually, in place of \eqref{interpol12}, we get the system
\begin{equation} \label{interpol22}
  \displaystyle
   \begin{cases}
     \frac{1-\rho}{p-2} = \frac{1}{2(q-1)},\\
      \frac{\rho}{p} + \frac{1-\rho}{3(p-2)}\,=\,\frac{1}{r(q-1)},
   \end{cases}
\end{equation}
whence one computes, exactly as before,
\begin{equation} \label{interpol25}
  K_\epsi=\frac{9+7\epsi}{9+3\epsi}.    
\end{equation} 
and, finally,
\begin{align} \no
  q & = \frac12 \frac{p-2}{1-\rho} + 1
    = \frac12 \Big(1 + \frac{p}{p-2} K\ee \Big) (p-2) + 1\\
  \label{interpol28}
  & = \frac{ K\ee + 1 }2 p = \frac{9+5\epsi}{9+3\epsi} p = H p,
\end{align}
which, exactly as before, is larger than one. Hence, the procedure 
continues as before, with small variations in the numerical values of 
the indices. Of course, the $L^{3+\epsi}$ regularity of 
the initial datum is used since we need to take $p=p_0=3+\epsi$ at the
first iteration (for smaller values of $p$ we get no summability gain
from the gradient term).
\end{proof}
\noindent%
{\bf Conclusion of proof.}~~%
Thanks to estimate~\eqref{st-51}, $\teta$ satisfies the first \eqref{inpiu-2}
for any $S\ge 1$ (where $M$ is the quantity on the \rhs\ of \eqref{st-51}
which is independent of $S$). Analogously, we have $\teta(S)\in L^{3+\epsi}$ for 
(almost) all $S\ge 1$. Moreover, combining \eqref{st-energy} and \eqref{st-31b},
we have the second of \eqref{inpiu-2}, still with $F$ independent of $S$. We
then conclude applying the above Lemma over the generic interval $(S,S+2)$, with
$S\ge 1$.
\beos\label{regofinale}
 Of course, with \eqref{asym-u} and \eqref{asym-teta} at our disposal,
 equation~\eqref{calore} is both nonsingular and nondegenerate. Consequently,
 we can prove, with standard tools, further regularization properties
 of solutions. In the non-conserved case, thanks to \eqref{separ}, 
 also the (possibly) singular character of $b$ is lost. Thus, the smoothness
 of the solution is limited only by the differentiability properties
 of $b$. For instance, if $b\in C^\infty$, then also the solution is 
 infinitely differentiable for strictly positive times.
\eddos

%%%%%%%%%%%%%%%%%%%%%%%%%%%%% References %%%%%%%%%%%%%%%%%%%%%%%%%%    

%\section{References}

%\vspace{1cm}

\end{document}